\documentclass{amsart}

\usepackage{amssymb}

\newtheorem{thm}{Theorem}

\newtheorem{lem}[thm]{Lemma}

\newtheorem{prop}[thm]{Proposition}

\theoremstyle{definition}
\newtheorem{defn}[thm]{Definition}

\newtheorem{exmp}[thm]{Example}

\newtheorem{ques}[thm]{Question}    

\newtheorem*{ack}{Acknowledgments}      

\newtheorem{defn-thm}[thm]{Definition--Theorem}  
\newtheorem{defn-lem}[thm]{Definition--Lemma}  

\theoremstyle{remark}

\let \cedilla =\c
\renewcommand{\c}[0]{{\mathbb C}}

\renewcommand{\o}[0]{{\mathcal O}}
\newcommand{\z}[0]{{\mathbb Z}}

\renewcommand{\r}[0]{{\mathbb R}}

\newcommand{\p}[0]{{\mathbb P}}

\newcommand{\q}[0]{{\mathbb Q}}
\newcommand{\map}[0]{\dasharrow}
\newcommand{\qtq}[1]{\quad\mbox{#1}\quad}

\newcommand{\tsum}[0]{\textstyle{\sum}}

\begin{document}
\bibliographystyle{amsalpha}

\hfill\today

\title[Continuous  rational functions]{Continuous  rational functions on \\
real and $p$-adic varieties}
\author{J\'anos Koll\'ar and Krzysztof Nowak}

\maketitle

A {\it rational function} on $\r^n$ is a
quotient of two polynomials
$$
f(x_1,\dots, x_n):=\frac{p(x_1,\dots, x_n)}{q(x_1,\dots, x_n)}.
$$

Strictly speaking, a  rational function $f$
is not really a function on $\r^n$ in general
since it is defined only on the dense open set
where $q\neq 0$. Nonetheless, even if $q$ vanishes at some
points of $\r^n$, it can happen that there is an
everywhere defined continuous  function $f^c$
that agrees with $f$ at all points where $f$ is defined.
Such an $f^c$ is unique. For this reason it is customary
to  identify $f$ with $f^c$
and call  $f$ itself  a {\it continuous rational function} on $\r^n$.
For instance,
$$
\frac{p(x_1,\dots, x_n)}{x_1^{2m}+\cdots+x_n^{2m}}
$$
is  a continuous rational function on $\r^n$ if every monomial
in $p$ has degree $>2m$.

The above definition of  continuous rational functions makes sense
on any real algebraic variety $X$, as long as
the  open set where $f$ is defined is dense in $X$ in the  Euclidean topology.
This condition always holds on smooth varieties,
and, more generally, if the singular set is contained in the
  Euclidean closure of the smooth locus.

The aim of this note is to consider four
basic problems on  continuous rational functions.

\begin{ques}\label{question12}
Let $X$ be a real algebraic variety and $Z$ a closed
 subvariety.
\begin{enumerate}
\item
Let $f$ be a  continuous  rational function defined on $X$.
Is the restriction  $f|_Z$ a   rational function on $Z$?
\item  Let $g$ be  a continuous  rational function defined  on $Z$.
Can one extend it
 to a  continuous  rational function on $X$?
\item Assume that $X\setminus Z$ is  Zariski dense in $X$.
Is $f$ uniquely determined by  the restriction $f|_{X\setminus Z}$?
\item Which systems of linear equations
$$
\tsum_j f_{ij}\cdot y_j = g_i \qquad i=1,\dots, m
$$
have continuous  rational solutions
where the $g_i, f_{ij}$ are polynomials (or rational functions) on $X$.
\end{enumerate}
\end{ques}

Section 1 contains a series of counterexamples to Questions
\ref{question12}.1--3. Proposition~\ref{rtl.restrict.rtl} shows
that the answer to Question \ref{question12}.1 is always yes if
$X$ is smooth. In Section~2 this leads to the
 introduction of   the notion of
{\it  hereditarily rational functions.}
These have the good properties that one would expect
based on the smooth case. In Theorem \ref{crat.extend.thm} this
concept is used
 to give a
complete answer to Question \ref{question12}.2.

The key result   says that a continuous
hereditarily rational function defined on a closed subvariety
extends to a continuous
hereditarily rational function on the ambient variety.
We give two proofs of this result.

The first version,  Proposition~\ref{crat.extend.thm.2},
is quite explicit and gives optimal control over the
regularity of the extension. It uses basic properties of semialgebraic sets,
thus it works only over $\r$.
The second variant, discussed in Section 3,
works for all locally compact fields.
However, it relies on a rather strong version of resolution of
singularities and it does not control the regularity of the extension.

It is quite likely that similar results hold for more general
complete valued fields. We hope to consider these in the future,
as well as an  answer to  Question \ref{question12}.4.

Topological
aspects of continuous rational maps between smooth algebraic varieties were
investigated in \cite{kuch}.
 Hereditarily rational functions are also studied in
\cite{fhmm}. Applications of hereditarily rational functions are given
in the  paper~\cite{kuch-kur} discussing
stratified-algebraic vector bundles.


\section{Examples}

The following example shows that the answer to
 Questions \ref{question12}.1--2 is not always positive.

\begin{exmp} \label{nonextend.exmp}
Consider the  surface $S$ and the rational function
$f$  given by
$$
S:=\bigl(x^3-(1+z^2)y^3=0)\subset \r^3
\qtq{and}  f(x,y,z):=\tfrac{x}{y}.
$$
We claim that
\begin{enumerate}
\item $S$ is a real analytic submanifold of $\r^3$,
\item $f$ is defined away from the $z$-axis,
\item $f$ extends to  a
real analytic function $f^c$ on $S$, yet
\item the restriction of $f^c$ to the $z$-axis is not rational and
\item $f$ can not be extended to a  continuous rational function on $\r^3$.
\end{enumerate}

Proof. Note that
$$
x^3-(1+z^2)y^3=\bigl(x-(1+z^2)^{1/3}y\bigr)
\bigl(x^2+(1+z^2)^{1/3}xy+(1+z^2)^{2/3}y^2\bigr).
$$
The first factor defines $S$ as a real analytic submanifold of $\r^3$.
The second factor only vanishes along the $z$-axis,
which is contained in $S$.
Therefore $x-y\sqrt[3]{1+z^2}$ vanishes on $S$, hence
$$
f^c|_S=\tfrac{x}{y}|_S=\sqrt[3]{1+z^2}|_S\qtq{and so} f(0,0,z)=\sqrt[3]{1+z^2}.
$$

Assume finally that $F$ is a continuous rational function on
$\r^3$ whose restriction to $S$ is $f$.
Then $F$ and $f$ have the same restrictions to the $z$-axis.
We show below that the  restriction of any
continuous rational function $F$ defined on $\r^3$
to the $z$-axis is a rational function.
Thus $F|_{z-{\rm axis}}$
 does not equal $\sqrt[3]{1+z^2}$, a contradiction.

To see the claim, write
$F={p(x, y, z)}/{q(x, y, z)}$
where $p,q$ are polynomials. We may assume that they are relatively prime.
Since $F$ is continuous everywhere, $x$ can not divide $q$.
Hence
$F|_{(x=0)}={p(0, y, z)}/{q(0, y, z)}$.
By canceling common factors, we can write this as
$F|_{(x=0)}={p_1(y, z)}/{q_1( y, z)}$
where $p_1,q_1$ are relatively prime polynomials.
As before, $y$  can not divide $q_1$, hence
$F|_{z-{\rm axis}}={p_1(0, z)}/{q_1(0, z)}$
is a rational function.

(Note that we seemingly have not used the continuity of $F$:
for any rational function $f(x,y,z)$
the above procedure defines a rational function
on the ${z-{\rm axis}}$. However, if we use
$x,y$ in reverse order, we could get a
different   rational function. This happens, for instance,
for $f=x^2z/(x^2+y^2)$.
Here  $\bigl(f|_{(x=0)}\bigr)|_{z-{\rm axis}}=0$ and
$\bigl(f|_{(y=0)}\bigr)|_{z-{\rm axis}}=z$.)
\end{exmp}
\medskip

In the above example, the problems arise since $S$ is not normal.
However, the key properties (\ref{nonextend.exmp}.4--5) can also be realized
on a normal hypersurface.

\begin{exmp} \label{nonextend.exmp.2} Consider
$$
X:=\bigl((x^3-(1+t^2)y^3)^2+z^6+y^7=0\bigr)\subset \r^4
\qtq{and}  f(x,y,z,t):=\tfrac{x}{y}.
$$
We easily see that the singular locus is the
$t$-axis.
$X$ is normal since a hypersurface in a smooth variety is normal
iff its singular set has codimension $\geq 2$; see for instance
\cite[Chap.III.C,Prop.9]{ser-al} or
\cite[Thm.11.2]{eis-ca}.

Let us blow  up the $t$-axis. There is one relevant chart,
where $x_1=x/y, y_1=y, z_1=z/y$. We get the smooth 3-fold
$$
X':=\bigl((x_1^3-(1+t^2))^2+z_1^6+y_1=0\bigr)\subset \r^4.
$$
Each point $(0,0,0,t)$ has only 1 preimage in $X'$,
given by  $\bigl(\sqrt[3]{1+t^2}, 0,0,t\bigr)$
and the projection $\pi:X'\to X$ is a homeomorphism.
Thus  $X$ is a topological  manifold, but it is not a
differentiable submanifold of $\r^4$.

Since $f\circ \pi=x_1$ is a regular function,
we conclude that $f\circ \pi$ extends to a continuous
(even regular) function $(f\circ \pi)^c$ on $X'$.
Since $\pi$ is a homeomorphism,
we get the continuous function
$f^c:=(f\circ \pi)^c\circ \pi^{-1}$
 on $X$ extending $f$.
By construction,   $f(0,0,0,t)=\sqrt[3]{1+t^2}$,
thus, as before, $f$ can not be extended to a
continuous rational function on
$\r^4$.

For any $m\geq 1$ we get  similar examples
of normal hypersurfaces and rational functions
$$
X_m:=\bigl((x^3-(1+t^2)y^3)^2+z_1^6+\cdots+z_m^6+y^7=0\bigr)\subset \r^{3+m}
\qtq{and}  f:=\tfrac{x}{y}.
$$
Note that for all $m$, the singular set is still the $t$-axis
and, for $m\gg 1$,
 the $X_m$  have rational, even terminal singularities
(see \cite{km-book} for the definitions of these singularities).
In fact, we do not know any natural class of singularities
(other than smooth points) where  Questions \ref{question12}.1--2
have a positive answer.
\end{exmp}

In order to elucidate  Question \ref{question12}.3,
next we give an example of a
continuous  rational function $f$ on $\r^3$
and of an irreducible algebraic surface $S\subset \r^3$
such that $f|_S$ is zero on a Zariski dense, Zariski  open subset of $S$
yet $f^c|_S$ is not identically zero.

\begin{exmp} \label{tail.rest.exmp}
Consider the rational function
$$
f(x,y,z):=z^2\cdot \frac{x^2+y^2z^2-y^3}{x^2+y^2z^2+y^4}.
$$
Its only possible discontinuities
are along the $z$-axis. To analyze its behavior there, rewrite it as
$$
f=z^2-y(1+y)\cdot \frac{y^2z^2}{x^2+y^2z^2+y^4}.
$$
The fraction is  bounded by $1$, hence
$f$ extends to a  continuous  function $f^c$
and $f^c(0,0,z)=z^2$.

Our example is the restriction of $f$ to  the surface
$S:=\bigl(x^2+y^2z^2-y^3=0\bigr)\subset \r^3$.
The image of the projection of $S$ to the $(y,z)$-plane is
 the inside of the parabola
$(y\geq z^2)\subset \r^2_{yz}$.
Topologically, $S$ has 2 parts. One is the $z$-axis, which is also the
singular locus of $S$, and
the other part $S^*$ is the Euclidean closure of the smooth locus of $S$.
The 2 parts intersect only at the origin.
Thus $S^*$ is Zariski dense but not
Euclidean dense in $S$.

We see that $f^c$ vanishes on $S^*$ but not on the $z$-axis.

More generally, let $g(z)$ be any rational function without real poles.
Then $g(z)f(x,y,z)$ vanishes on $S^*$
and its restriction to the $z$-axis is $z^2g(z)$.

(The best known example of a
surface with a Zariski dense open set which is not  Euclidean  dense
is the
Whitney umbrella  $W:=(x^2=y^2z)\subset \r^3$.
We can take $W_1:=(x=y=0)$ and $W_2:=W$.
The  Euclidean closure of $W_2\setminus W_1$  does not contain the ``handle''
$(x=y=0, z<0)$.
In this case, a
 continuous  rational function is determined by
its   restriction to $W_2\setminus W_1$.
The Euclidean closure of $W_2\setminus W_1$
contains the half line $(x=y=0, z\geq 0)$, and
a rational function on a line is determined by
its restriction to any interval.)
\end{exmp}

The next example shows the two natural ways of pulling back
continuous rational functions by a  morphism
can be  different.

\begin{exmp}[Two pull-backs]\label{two.pb.exmp}
Let $\pi:X'\to X$ be a  morphism of
real varieties.
Assume for simplicity that $X,X'$ are irreducible and $\pi$ is birational.
If $f$ is a continuous rational function on
$X$, then one can think of the
pull-back of $f$ or of $f^c$ to $X'$
in at least two different ways.

First, $f^c\circ \pi$ is the composite of
two continuous maps, hence it is a
continuous function.
Second, one can view $f\circ \pi$ as a rational function
on $X'$. Let $U'\subset X'$ be the open set where $f\circ \pi$
is regular. Let us denote the resulting regular function
by $(f\circ \pi)^r:U'\to \r$.

The following example shows that
$(f\circ \pi)^r$ and $\bigl(f^c\circ \pi\bigr)|_{U'}$
 can be different.

Take $X=\r^2$ with $f=x^3/(x^2+y^2)$.
Note that $f^c(0,0)=0$.
Blow up $(x^3,x^2+y^2)$ to obtain
$X'\subset \r^2\times \r\p^1$. The first projection
$\pi:X'\to X$ is an isomorphism away from the origin
and $\pi^{-1}(0,0)\cong \r\p^1$.

The first interpretation above gives a continuous function
$f^c\circ \pi$ which vanishes along $\pi^{-1}(0,0)$.

The second  interpretation views $f\circ \pi$
as a rational map
$$
f\circ \pi: X'\setminus \pi^{-1}(0,0) \map \r\p^1
$$ which is in fact
regular on $X'$.
Thus it extends to a continuous (even regular) map
$(f\circ \pi)^c$ on $X'$ whose restriction to
 $\pi^{-1}(0,0)$ is an isomorphism
$\pi^{-1}(0,0)\cong \r\p^1$.

This confusion is possible only because
$X'_0:=\pi^{-1}\bigl(\r^2\setminus{(0,0)}\bigr)$ is not Euclidean dense in $X'$.
Its  Euclidean  closure contains only one point of $\pi^{-1}(0,0)\cong \r\p^1$.
The two versions of $f\circ \pi$ agree on $X'_0$, hence also
on its Euclidean closure, but not everywhere.

In general, let $\pi:X'\to X$ be a  morphism of
real varieties and $f$ a rational function on $X$ that is regular
on an open set $U\subset X$. Let $U'\subset X'$ be the Euclidean closure
of $\pi^{-1}(U)$; it  is a semialgebraic subset of $X'$.
All possible definitions of a  continuous pull-back of
$f$   agree on $U'$ but, as the above example shows,
they may be different on $X'\setminus U'$.

\end{exmp}

Finally  we turn to  Question \ref{question12}.4
for a single equation
$$
\tsum_i\ f_i({\mathbf x})\cdot y_i= g({\mathbf x}),
$$
where $g$ and the $f_i$ are polynomials in the variables
${\mathbf x}=(x_1,\dots, x_n)$.
Such equations have a solution where the $y_i$ are
rational functions provided not all of the $f_i$ are
identically zero. The existence of a solution where
 the $y_i$ are
continuous functions is studied in \cite{fef-kol}
and \cite{k-hoch}.

One could then hope to prove that if there is a
continuous  solution then there is also
a continuous rational solution.
 \cite[Sec.2]{fef-kol} proves that  if there is a
continuous  solution then
there is also a  continuous semialgebraic solution.
The next example shows, however, that in general
there is no  continuous rational solution.

\begin{exmp}\label{main.exmp}
We claim that the linear equation
$$
 x_1^3x_2\cdot y_1+ \bigl(x_1^3-(1+x_3^2)x_2^3\bigr)\cdot y_2= x_1^4
\eqno{(\ref{main.exmp}.1)}
$$
has a continuous semialgebraic solution but no
continuous rational solution.

A  continuous semialgebraic solution is given by
$$
y_1=(1+x_3^2)^{1/3}\qtq{and}
y_2=\frac{x_1^3}{x_1^2+(1+x_3^2)^{1/3}x_1x_2+(1+x_3^2)^{2/3}x_2^2}.
\eqno{(\ref{main.exmp}.2)}
$$
(Note that
$x_1^2+(1+x_3^2)^{1/3}x_1x_2+(1+x_3^2)^{2/3}x_2^2\geq \frac12 (x_1^2+x_2^2)$,
so $|y_2|\leq 2x_1$ and it is indeed continuous.)
A solution by rational functions is
$y_1=x_1/x_2$ and $y_2=0$.

To see that (\ref{main.exmp}.1)
has no continuous rational solution,
restrict any solution $(y_1, y_2)$ to the semialgebraic surface
$S:=(x_1-(1+x_3^2)^{1/3}x_2=0)$.
Since $x_1^3-(1+x_3^2)x_2^3$ is identically zero on $S$,
 we conclude that $y_1|_S=x_1^4/(x_1^3x_2)|_S=(x_1/x_2)|_S$.
The latter is equal to $\sqrt[3]{1+x_3^2}$, thus
$$
y_1|_{x_3-{\rm axis}}=\sqrt[3]{1+x_3^2},
$$
which is not a rational function.
As we saw in Example \ref{nonextend.exmp},
this implies that
$y_1$ is not a rational function.
\end{exmp}

 \section{Hereditarily rational functions}

First we show that  Question \ref{question12}.1 has a positive answer
 on  smooth varieties.

\begin{prop}\label{rtl.restrict.rtl}
Let $X$ be a real algebraic variety and $Z$ an irreducible
subvariety that is not contained in the singular locus of $X$.
Let $f$ be a   rational function on $X$ with continuous extension $f^c$.
Then there is a Zariski dense open subset $Z^0\subset Z$ such that
$f^c|_{Z^0}$ is a   regular function.
\end{prop}

Proof. By replacing $X$ with a suitable open subvariety, we may assume that
$X$ and $Z$ are both smooth.

Assume first that $Z$  has codimension 1.
Then the local ring $\o_{X,Z}$ is a principal ideal domain
\cite[Sec.II.3.1]{shaf};
let $t\in \o_{X,Z}$ be a defining equation of $Z$.
We can write $f=t^mu$ where $m\in \z$ and $u\in \o_{X,Z}$ is a unit.
Here  $m\geq 0$ since $f$ does not have a pole along $Z$, hence $f$ is regular
 along  a Zariski dense open subset $Z^0\subset Z$.
Thus $f|_{Z^0}$ is a  regular  function.

If $Z$ has higher codimension,
 note that $Z$ is a local complete intersection
at its smooth points \cite[Sec.II.3.2]{shaf}. That is,  there is a
sequence of subvarieties
$Z_0\subset Z_1\subset \cdots\subset Z_m=X_0\subset X$
where each $Z_i$ is a smooth hypersurface in $Z_{i+1}$ for $i=0,\dots, m-1$
and $Z_0$ (resp.\ $X_0$) is open and dense in $Z$ (resp.\ $X$).
We can thus restrict $f^c$ to $Z_{m-1}$, then to $Z_{m-2}$ and so on,
until we get that $f^c|_{Z_0}$ is regular.
(As we noted in Example \ref{main.exmp},
 for any rational function $f$ he above procedure defines
 a regular function
$f|_{Z_0}$, but it depends on the choice of the chain
$Z_1\subset  \cdots\subset Z_m$.)
\qed

\medskip

This suggests that we should focus on those functions
for which  Question \ref{question12}.1 has a positive answer.
Then we show that for such functions
Question \ref{question12}.2 also has a positive answer.

\begin{defn} \label{her.rtl.defn}
Let $X$ be a real algebraic variety and  $h$
 a  continuous  function on $X$.
We say that $h$ is {\it hereditarily rational}
if  every irreducible, real subvariety $Z\subset X$
has a Zariski dense open subvariety $Z^0\subset Z$
such that
$h|_{Z^0}$ is a  regular function on $Z^0$.
(A  function $h$ is {\it regular at} $x\in X$
if one can write $h=p/q$ where $p,q$ are polynomials and
$q(x)\neq 0$.  It is called {\it regular} if it is
regular at every point of $X$.)

Examples \ref{nonextend.exmp} and  \ref{nonextend.exmp.2}
show that not every  continuous  rational function
is  hereditarily rational.

If $h$ is rational, there is a Zariski dense open set $X^0\subset X$
such that  $h|_{X^0}$ is regular.
 If $h$ is  continuous and hereditarily rational,
we can repeat this process with the restriction of $h$ to
$X\setminus X^0$, and so on. Thus we conclude that
 a continuous function $h$ is  hereditarily rational
 iff there is a sequence of closed subvarieties
$\emptyset=X_{-1}\subset X_0\subset\cdots\subset X_m=X$
such that for  $i=0,\dots, m$ the restriction of $h$ to
 $ X_i\setminus X_{i-1}$ is regular.

 If it is convenient, we can also assume that each
$X_i\setminus X_{i-1}$ is smooth of pure dimension $i$.

 Proposition \ref{rtl.restrict.rtl} says  that on a smooth variety
every continuous  rational function  is
hereditarily rational.

The pull-back of a hereditarily rational function
by any morphism is again a
(continuous and) hereditarily rational function.
\end{defn}

The following result shows that
 hereditarily rational functions constitute the right class for
 Question \ref{question12}.2.

\begin{thm} \label{crat.extend.thm}
Let $Z$ be a real algebraic variety and
 $f$  a    rational function on $Z$ with continuous extension $f^c$.
The following are equivalent.
\begin{enumerate}
\item $f^c$ is hereditarily rational.
\item For every real algebraic variety $X$ that
contains $Z$ as a closed subvariety,
 $f^c$ extends to a   (continuous and)
hereditarily  rational function $F$ on $X$.
\item For every real algebraic variety $X$ that
contains $Z$ as a closed subvariety,
 $f$ extends to a continuous   rational function $F$ on $X$.
\item Let $X_0$ be a smooth real algebraic variety that
contains  $Z$ as a closed subvariety. Then
 $f$ extends to a continuous  rational function $F_0$ on $X_0$.
\end{enumerate}
\end{thm}

Proof. It is clear that (2) $\Rightarrow$ (3) $\Rightarrow$ (4)
and (4)  $\Rightarrow$ (1) holds by Proposition \ref{rtl.restrict.rtl}.
Thus we need to show that
(1)  $\Rightarrow$ (2).

 We can embed $X$ into a smooth  real algebraic variety $X'$.
If we can extend $f$ to $X'$ then its restriction to $X$ gives the
required extension. Thus we may assume to start with that $X$ is smooth
(or even that $X=\r^N$ for some $N$).

We  prove the following  more precise version.

\begin{prop} \label{crat.extend.thm.2}
Let $X$ be a smooth real algebraic variety and
$W\subset Z\subset X$ closed subvarieties.
Let $f$ be a continuous, hereditarily  rational function  on $Z$
that is regular on $Z\setminus W$.
Then $f$ extends to a   continuous,
hereditarily  rational function $F$ on $X$
that is regular on $X\setminus W$.
\end{prop}

Proof. We use induction on $\dim Z$. The case $\dim Z=0$ is obvious.

If $W$ is replaced by a smaller set, the assertion gets stronger.
Hence we may assume that $W\subset Z$ is the smallest
set such that $f$ is regular on $Z\setminus W$.
Since every rational function is regular on a Zariski dense open set,
 $W$ is nowhere dense in $Z$. In particular,
$\dim W<\dim Z$.

Since $f$ is hereditarily  rational, $f^c|_W$ is also
hereditarily  rational. Thus, by induction,
there is a hereditarily  rational function $F_1$ on $X$
that is regular on $X\setminus W$ and such that $F_1^c|_W=f^c|_W$.
Set $f_2:=f-F_1|_W$. Then $f_2^c$ vanishes on $W$ and it is
enough to show that  $f_2$ extends to a   continuous,
hereditarily  rational function $F_2$ on $X$
that is regular on $X\setminus W$.

Let $I_2$ be the ideal of those regular functions $\phi$ on $Z$
such that $\phi f_2$ is regular. Let $\phi_1,\dots, \phi_r$ be
a set of generators of $I_2$. Set $q:= \phi_1^2+\cdots + \phi_r^2$
then $p:=qf_2$ is regular, $q$ vanishes precisely
at the points where $f$ is not regular and $f_2=p/q$.
Note further that $W$ is the zero set of $q$ and the zero set of
$p$ contains $W$.

A regular function on a closed subvariety $Z$
of a real algebraic variety $X$ always extends to
 a regular function on the whole variety.
(This follows, for instance,  from  \cite[4.4.5]{bcr}
or see Lemma \ref{ext.reg.lem} for a more general argument.)

Thus $p$ (resp.\ $q$) extend to regular functions
$P$ (resp.\ $Q$) on $X$. Usually $P/Q$ is not continuous
but we will improve it in two steps.

Let $H$ be a regular function on $X$ whose zero set is $Z$.
Then $Q^2+H^2$ vanishes only on $W$ and its restriction to
$Z$ equals $q^2$. Thus
$$
G:=\frac{PQ}{Q^2+H^2}
$$
is a rational  function on $X$
that is regular on $X\setminus W$ and whose restriction
to $Z\setminus W$ equals $f_2$. However, usually
it is not continuous on $X$ near $W$ hence
we need one more correction.
\medskip

{\it Claim.} For $n\gg 1$ the function
$$
F_{2n}:=G\cdot \frac{Q^{2n}}{Q^{2n}+H^2}=
\frac{PQ}{Q^2+H^2}\cdot \frac{Q^{2n}}{Q^{2n}+H^2}
$$
is continuous on $X$.
\medskip

It is clear that the restriction of $F_{2n} $
to $Z\setminus W$ equals $f_2$, thus, once the Claim is proved,
$F:=F_1+F_{2n}$ satisfies all the requirements.

In order to prove the claim, we work on the variety
$\pi:X_1\to X$ obtained by blowing up the ideal $(PQ, Q^2+H^2)$.
Equivalently, $X_1$ is the Zariski closure of the graph of $G$
in  $X\times \r\p^1$.

Since $W$ is the common zero set  $(PQ=Q^2+H^2=0) $,
the exceptional set of $\pi$ is
$E:=\pi^{-1}(W)\cong W \times \r\p^1$
and $\pi$ is an isomorphism over $X\setminus W $.
We prove that
$$
\bigl(F_{2n}|_{X\setminus W}\bigr)\circ \pi
$$
extends to a continuous function $\bigl(F_{2n}\circ \pi\bigr)^c$
that vanishes on  the exceptional set $E$.
Thus $\bigl(F_{2n}\circ \pi\bigr)^c$ descends to a
continuous function on $X$ that vanishes on  $W$.

The existence of such a continuous extension is a local question
on $X_1$ and we use different arguments on different charts.
A slight complication is that one of our charts is only semialgebraic.

Let $Z^*\subset X_1$ be the Euclidean closure of
$\pi^{-1}(Z\setminus W)$. Note that  $Z^*$ is a closed semialgebraic subset
 but it is usually not real algebraic.

On $X_1$ one can identify the rational function $G\circ \pi$
with the  restriction of the second projection
$\pi_2:W \times \r\p^1\to \r\p^1$.

On $Z^*$ we thus have a continuous function $f_2^c\circ \pi$
that vanishes on  $Z^*\cap E$ and a regular (hence continuous) map
$\pi_2|_{Z^*}:Z^*\to  \r\p^1$. Using the natural identification
$\r^1= \r\p^1\setminus\{\infty\}$, these two agree on the open set
$\pi^{-1}(Z\setminus W)$, hence they agree on $Z^*$.
(Note that, as in Example \ref{two.pb.exmp},
 the two functions  might not agree on the Zariski closure of
$\pi^{-1}(Z\setminus W)$; this is why we work with $Z^*$.)
Thus $\pi_2^{-1}(\infty)$ is disjoint from
$Z^*$, and  hence there is
 a Zariski open  neighborhood $U^* \subset X_{1}$ of $Z^*$
such that $\pi_2$ defines a regular function on $U^*$.
This gives the  extension of
$\bigl(G|_{X\setminus W}\bigr)\circ \pi$
to a regular function on $U^*$ that vanishes along $Z^*\cap E$.
(It is not clear that  $(G \circ \pi)|_{U^*}$
vanishes along $U^*\cap E$.)

Note further that
 $$
\frac{Q^{2n}}{Q^{2n}+H^2}\circ \pi
$$
is a bounded regular function on $X_1\setminus E$.
Therefore the restriction of the product
$$
\bigl(F_{2n}\circ \pi\bigr)\bigr|_{U^*}=
\bigl(G\circ \pi \bigr)\cdot  \Bigl(\frac{Q^{2n}}{Q^{2n}+H^2}\circ \pi\Bigr)
\Bigr|_{U^*}
$$
extends to a  function that vanishes and
is continuous at every point of $Z^*\cap E$.
However,
 we have not proved so far that it is defined along $(U^*\cap E)\setminus Z^*$.

The other chart is $V^* := X_{1} \setminus Z^*$;
it is semialgebraic and Euclidean open in $X_1$.
Here we write
$F_{2n}\circ \pi$ in the form
$$
F_{2n}\circ \pi=
(P\circ \pi)\cdot \Bigl(\frac{Q^{2n-1}}{H^2}\circ \pi\Bigr)
\cdot \Bigl(\frac{Q^2}{Q^2+H^2}\circ \pi\Bigr)\cdot
\Bigl(\frac{H^2}{Q^{2n}+H^2}\circ \pi\Bigr).
$$
The last 2 factors are bounded regular functions on $V^*\setminus E$ and
$P$ is a  regular function that vanishes
along  $E\cap V^*$.

Note that on $V^*$ the function $H \circ \pi$ vanishes only along $E$
and  $Q \circ \pi$ also vanishes  along $E$.
We can thus apply Theorem \ref{lojas.thm}
to conclude that $\bigl(Q^{2n-1}/H^2\bigr)\circ \pi$
 extends to a  continuous (and semialgebraic)
function on $V^*$ for $n\gg 1$. Thus $ F_{2n}\circ \pi$
extends to a  continuous function on $V^*$
that vanishes along $E\cap V^*$
for $n\gg 1$.

Putting the two charts together we conclude that
 $ F_{2n}\circ \pi$
extends to a  continuous function on $X_1$ that vanishes along $E$.
\qed
\medskip

 We have used the following version
of the \L{}ojasiewicz inequality given in \cite[Thm.2.6.6]{bcr},
see also  \cite{Now0,Now}.

\begin{thm}[\L{}ojasiewicz Inequality]\label{lojas.thm}
Let $V$ be a locally closed, semialgebraic subset of $\r^N$ and
 $\phi,\psi: V \to \r$  continuous
semialgebraic functions.
Assume that
$ \{ x \in V: \phi(x)=0 \} \subset \{ x \in V: \psi(x)=0 \}$.
Then there exist a positive integer $n$ and a continuous semialgebraic function
$\rho$ such that $\psi^n = \rho\phi$.\qed
\end{thm}

The following example shows that, even in very simple cases,
the extension of continuous rational functions is not entirely trivial.

\begin{exmp} Let $X=\r^2$ and $Z\subset \r^2$ the cuspidal cubic
with equation $(x^2-y^3=0)$. Consider the rational function
$f(x,y):=y^2/x$.
$Z$ can be parametrized as $x=t^3, y=t^2$ and then
$f(t^3, t^2)=t$ is clearly continuous.

First we claim that there are no regular functions
$P,Q$ such that $P|_Z=y^2, Q|_Z=x$ and $P/Q$ is continuous.
In fact, $P/Q$ can not even be bounded.
Indeed, any such extension would be of the form
$$
P=y^2+P_1(x^2-y^3)\qtq{and} Q=x+Q_1(x^2-y^3)
$$
where $P_1, Q_1$ are regular. Thus
$$
\frac{P}{Q}\Bigr|_{y-\mbox{axis}}=\frac{y^2-y^3P_1}{-y^3Q_1}=
\frac{1-yP_1}{-yQ_1}
$$
has a pole at $y=0$.

Our first two improvements are
$$
\frac{y^2x}{x^2+(x^2-y^3)^2}\qtq{and}
\frac{y^2x}{x^2+(x^2-y^3)^2}\cdot \frac{x^2}{x^2+(x^2-y^3)^2}.
$$
For both of these, the limit along the curve $(t^2, t)$ is 1,
hence they are not  continuous at the origin.
The next improvement is
$$
F_4:=\frac{y^2x}{x^2+(x^2-y^3)^2}\cdot \frac{x^4}{x^4+(x^2-y^3)^2}.
$$
This turns out to be continuous as can be seen either by two blow-ups
or by the following direct computation.
Performing the change of variables
$x=u^{3}, \ \ y=u^{2}v$ we get
$$
\bigl|F_4(u^3, u^2v)\bigr|=
\Bigl|\frac{uv^2}{(1+u^6(1-v^3)^2) \cdot (1 + (1-v^3)^2)}\Bigr|
\leq |u|\cdot \frac{v^2}{1 + (1-v^3)^2}.
$$
The last factor is uniformly bounded by a constant $C$ for $v\in \r$,
thus we conclude that
$$
\bigl| F_4(x,y)\bigr| \leq C |x|^{1/3}.
$$
\end{exmp}


\section{Varieties over $p$-adic fields}

Let $K$ be any topological field.
The $K$-points $X(K)$ of any $K$-variety inherit from $K$
a topology, called the $K$-topology. One can then consider
rational functions $f$ on $X$ that are continuous on $X(K)$.
This does not seem to be a very interesting notion in general,
unless $K$ satisfies the following {\it density property.}
\begin{enumerate}
\item[(DP)] If $X$ is smooth, irreducible
 and $\emptyset\neq U\subset X$ is  Zariski open
then $U(K)$ is dense in $X(K)$ in the $K$-topology.
\end{enumerate}
Note that if (DP) holds for all smooth curves then it holds for
all varieties. Examples of such fields are  complete real
valued fields with non-discrete topology; this can be deduced by
means of Hensel's lemma for restricted formal power series over
real valued fields (see e.g.~\cite[Chap.III]{Bour}).

If $K$ is algebraically closed, for instance $K=\c$,
then every continuous rational function on a normal variety
is regular, so we do not get a new notion. In general, the study of
continuous rational functions leads to
the concept of {\it seminormality} and {\it seminormalization}; see
\cite{MR0239118, MR0266923} or \cite[Sec.10.2]{kk-singbook} for a recent
treatment.


The proofs of Section 2, except that of
Proposition~\ref{crat.extend.thm.2}, all work over such general
topological fields $K$. As for Proposition~\ref{crat.extend.thm.2},
there seem to be 3 issues.

$\bullet$ We need to show that  every rational
function that is regular at the $K$-points of a subvariety extends to
a rational
function that is regular at the $K$-points of the ambient  variety.
This is proved by a slight modification of the usual arguments that apply
when the base field is algebraically closed \cite[Sec.I.3.2]{shaf} or real
closed \cite[3.2.3]{bcr}.

$\bullet$ The proof of Proposition~\ref{crat.extend.thm.2} relied
on semialgebraic sets. It is doubtful that definable sets behave
sufficiently nicely for  arbitrary
topological fields, though probably one could use them for
$p$-adic fields. We go around this problem by using a
strong form of resolution of singularities and
transformation to a simple normal crossing; see \cite[Chap.III]{k-resol}  for
references and relatively short proofs.

$\bullet$ At the end of the proof we need to show that if
$\sigma:Y\to X$ is a birational morphism between smooth varieties
and $f_Y$ is a rational function that is continuous on $Y(K)$ and
constant on the fibers of $Y(K)\to X(K)$ then $f_Y$ descends to a
rational function $f_X$ that is continuous on $X(K)$. If $K$ is
locally compact then $Y(K)\to X(K)$ is proper, hence the
continuity of $f_X$ is clear. The similar assertion does not hold
for arbitrary fields satisfying (DP), even if $X$ is normal, but
we do not know what happens in the smooth case.

All the above steps work if $K$ is locally compact. These fields are
$\r$,  the $p$-adic fields $\q_p$ and their finite
extensions.

\begin{lem}[Extending regular functions]\label{ext.reg.lem}
Let $k$ be a field, $X$ an affine $k$-variety
and $Z\subset X$ a closed subvariety. Let $f$ be a rational function on
$Z$ that is regular at all points of $Z(k)$.
Then there is a  rational function $F$ on
$X$ that is regular at all points of $X(k)$ and such that
$F|_Z=f$.
\end{lem}

Proof. If $k$ is algebraically closed, then $f$ is regular on
$Z$ hence it extends to a regular function on $X$.

If $k$ is not algebraically closed, then, as an auxiliary step,
 we claim that there are
polynomials $G_{r}(x_1,\dots, x_r)$ in any number of variables
whose only zero on $k^r$ is $(0,\dots,0)$.



Indeed, take a polynomial
$$ g(t) = t^{d} + a_{1}t^{d-1} + \cdots + a_{d} \in k[t], \ \ d>1, $$
which has no roots. Then its homogenization
$$ G_{2}(x_{1},x_{2}) = x_{1}^{d} + a_{1}x_{1}^{d-1}x_{2} + \cdots + a_{d}x_{2}^{d} $$
is a polynomial in two variables we are looking for. Further, we
can recursively define polynomials $G_{r}$ by putting
$$ G_{r+1}(x_{1},\ldots,x_{r},x_{r+1}) :=
   G_{2}(G_{r}(x_{1},\ldots,x_{r}),x_{r+1}). $$

Now we construct the extension of $f$ as follows.

For every $z\in Z(k)$ we can write $f=p_z/q_z$ where $q_z(z)\neq
0$. After multiplying both $p_z,q_z$ with a suitable polynomial,
we can assume that $p_z, q_z$ are regular on $Z$ and then
 extend them to  regular functions on $X$.
By assumption, $\bigcap_{z\in Z(k)} (q_z=0)$ is disjoint from
$Z(k)$. Choose finitely many $z_1,\dots, z_m\in Z(k)$ such that
$$
\textstyle{\bigcap_{i=1}^m} \bigl(q_{z_i}=0\bigr)= \bigcap_{z\in
Z(k)} \bigl(q_z=0\bigr).
$$
Let $q_{m+1},\dots, q_r$ be defining equations of $Z\subset X$.
Set $q_i:=q_{z_i}$ and $p_i:=p_{z_i}$ for $i=1,\dots, m$ and
$p_i=q_i$ for $i=m+1,\dots, r$. Write (non-uniquely)  $G_{r}=\sum
G_{ri}x_i$ and finally set
$$
F:=\frac{\sum_{i=1}^r G_{ri}(q_1,\dots, q_r)p_i}{G_{r}(q_1,\dots,
q_r)}.
$$
Since the $q_i$ have no common zero on $X(k)$, we see that $F$ is
regular at all points of $X(k)$.

Along $Z$, $p_i=fq_i$ for  $i=1,\dots, m$ by construction and
for  $i=m+1,\dots, r$ since then both sides are 0. Thus
$$
F|_Z:=\frac{\sum_{i=1}^r G_{ri}(q_1,\dots,
q_r)fq_i}{G_{r}(q_1,\dots, q_r)}|_Z= f\cdot \frac{\sum_{i=1}^r
G_{ri}(q_1,\dots, q_r)q_i}{G_{r}(q_1,\dots, q_r)}|_Z=f. \qed
$$

>From now on $K$ denotes a locally compact field
of characteristic 0 that is not algebraicaly closed.
That is, $K$ is either $\r$ or a  finite
extension of a  $p$-adic field $\q_p$.

\begin{prop} \label{crat.extend.thm.2.gen}
Let $X$ be a smooth  $K$-variety and $Z\subset X$ a closed
subvariety. Let $f$ be a continuous, hereditarily  rational
function  on $Z(K)$. Then $f$ extends to a   continuous,
hereditarily  rational function $F$ on $X(K)$.
\end{prop}

Proof. Let $W\subset Z$ denote the closed subvariety
consisting of those points $z\in Z$ such that either $f$ is not
regular at $z$ or $Z$ is singular at $z$.
As in the proof of Proposition~\ref{crat.extend.thm.2},
an induction argument reduces the problem to the
following special case.
\medskip

{\it Claim.}  The conclusion of
 Proposition~\ref{crat.extend.thm.2.gen} holds if we assume in addition that
 $f$ vanishes on $W(K)$.
\medskip

In order to prove the Claim, take an embedded resolution of singularities
$\sigma: X'\to X$ with the following properties.
\begin{enumerate}
\item $\sigma$ is an isomorphism over $X\setminus W$,
\item the birational transform $Z':=\sigma^{-1}_*(Z)$ is smooth,
\item the rational map  $f\circ \sigma$ restricts to a
morphism $f':Z'\to \p^1$,
\item the exceptional set $E:=\sigma^{-1}(W)$
is a simple normal crossing divisor and
\item $Z'\cup E$ is a simple normal crossing
subvariety.
\end{enumerate}
The last assertion means that for every point
$p\in Z'\cup E$ one can choose
local coordinates  $y_1,\dots, y_N$ such that, in a neighborhood of $p$,
$Z'=(y_1=\cdots=y_n=0)$ and the irreducible components
of $E$ are given by the equations $\{y_j=0 : j\in J\}$
for a suitable subset $J\subset \{n+1,\dots, N\}$.

Since $Z'$ is smooth, $f'|_{Z'(K)}$ agrees with the topological
pull-back  $f^c\circ \sigma$. Thus $f'|_{Z'(K)}$ is a regular
function on $Z'(K)$ that vanishes on $Z'(K) \cap E(K)$.

Consider next the rational function $g$ on $Z'\cup E$
whose restriction to $Z'$ equals $f'$ and to $E$
equals $0$. We claim that $g$ is regular at the $K$-points
of $Z'\cup E$. This is a local problem and the
only points in question are in $Z'\cap E$.

So pick a $K$-point $p\in Z'\cap E$, an open neighborhood $p\in
U\subset X'$  and local coordinates $y_1,\dots, y_N$ as above. We
know that $f'$ is regular on $Z'\cap U$, hence it extends to a
regular function $\bar f$ on $U$. Since $f'$ vanishes along
$Z'\cap E$, the extension  $\bar f$ also  vanishes along $Z'\cap
E\cap U$. Note that the ideal of $Z'\cap E\cap U$ is generated by
the functions  $y_1, \cdots, y_n, \prod_{j\in J}y_j$, hence we can
write  (non-uniquely)
$$
\bar f= p_1y_1+\cdots+ p_ny_n+ q\textstyle{\prod}_{j\in J}y_j
$$
where the $p_i$ and $q$ are regular on $U$.
Thus
$$
\bar f - p_1y_1-\cdots- p_ny_n = q\textstyle{\prod}_{j\in J}y_j
$$
is a regular
function near $p$ which vanishes on $E$ and whose restriction
to $Z'$ agrees with $f'$. Thus the restriction of
$\bar f$ to $Z'\cup E$ is regular and equals $g$.

Since $g$ is regular at all $K$-points of $Z'\cup E$,
by Lemma~\ref{ext.reg.lem} it extends to a rational function
$G$ on $X'$ that is regular on $X'(K)$.
Furthermore $G$ vanishes on the exceptional set $E$.

Finally, $\sigma: X'(K)\to X(K)$ is a proper surjection and $G$ is
constant on every fiber of  $\sigma|_{X'(K)}$. Thus $G$ descends
to the required continuous function $F$ on $X(K)$. \qed


\medskip

\begin{ack} We thank Ch.~Fefferman, M.~Jarden, F.~Mangolte and B.~Totaro
 for useful  comments and questions. Partial financial support  for JK was provided by  the NSF under
grant number DMS-0758275.  For KN, research was partially
supported by the Polish Ministry of Science and Higher Education
under grant number N N201 372336.
\end{ack}

\providecommand{\bysame}{\leavevmode\hbox to3em{\hrulefill}\thinspace}
\providecommand{\MR}{\relax\ifhmode\unskip\space\fi MR }
\providecommand{\MRhref}[2]{%
  \href{http://www.ams.org/mathscinet-getitem?mr=#1}{#2}
}
\providecommand{\href}[2]{#2}

\vskip0.3cm

\noindent Princeton University, Princeton NJ 08544-1000

\begin{verbatim}kollar@math.princeton.edu\end{verbatim}
\medskip

\noindent Institute of Mathematics, Jagiellonian University,

\noindent ul.~Profesora \L{}ojasiewicza 6, 30-348 Krak\'{o}w, Poland

 \begin{verbatim}nowak@im.uj.edu.pl\end{verbatim}

\end{document}